\documentclass[12pt,reqno]{amsart}
\usepackage{amssymb,verbatim,enumerate,ifthen}
\usepackage[mathscr]{eucal}
\def\N{\mathbb{N}}
\def\R{\mathbb{R}}
\def\Q{\mathbb{Q}}
\def\Z{\mathbb{Z}}

\def\A{\mathscr{A}}

\newtheorem{theorem}{Theorem}
\newtheorem*{theorem*}{Theorem}
\def\Thm#1#2{\ifthenelse{\equal{#1}{*}}{\begin{theorem*}#2\end{theorem*}}
             {\begin{theorem}\label{T#1}#2\end{theorem}}}
\newtheorem{Atheorem}{Theorem}

\def\THM#1#2{\begin{Atheorem}\label{T#1}#2\end{Atheorem}}
\def\thm#1{Theorem~\ref{T#1}}

\newtheorem{proposition}[theorem]{Proposition}
\newtheorem*{proposition*}{Proposition}
\def\Prp#1#2{\ifthenelse{\equal{#1}{*}}{\begin{proposition*}#2\end{proposition*}}
             {\begin{proposition}\label{P#1}#2\end{proposition}}}

\newtheorem{corollary}[theorem]{Corollary}
\newtheorem*{corollary*}{Corollary}
\def\Cor#1#2{\ifthenelse{\equal{#1}{*}}{\begin{corollary*}#2\end{corollary*}}
             {\begin{corollary}\label{C#1}#2\end{corollary}}}
\def\cor#1{Corollary~\ref{C#1}}

\newtheorem{lemma}[theorem]{Lemma}
\newtheorem*{lemma*}{Lemma}
\def\Lem#1#2{\ifthenelse{\equal{#1}{*}}{\begin{lemma*}#2\end{lemma*}}
             {\begin{lemma}\label{L#1}#2\end{lemma}}}
\def\lem#1{Lemma~\ref{L#1}}
\newtheorem{Alemma}{Lemma}

\def\LEM#1#2{\begin{Alemma}\label{L#1}#2\end{Alemma}}

\theoremstyle{definition}
\newtheorem{remark}[theorem]{Remark}
\newtheorem*{remark*}{Remark}
\def\Rem#1#2{\ifthenelse{\equal{#1}{*}}{\begin{remark*}\rm #2\end{remark*}}
             {\begin{remark}\label{R#1}\rm #2\end{remark}}}

\newtheorem{example}[theorem]{Example}
\newtheorem*{example*}{Example}
\def\Exa#1#2{\ifthenelse{\equal{#1}{*}}{\begin{example*}\rm #2\end{example*}}
             {\begin{example}\label{Ex#1}\rm #2\end{example}}}

\def\eq#1{{\rm(\ref{E#1})}}
\def\Eq#1#2{\ifthenelse{\equal{#1}{*}}
  {\begin{equation*}\begin{aligned}#2\end{aligned}\end{equation*}}
  {\begin{equation}\begin{aligned}\label{E#1}#2\end{aligned}\end{equation}}}

\def\dom{\mathop{\hbox{\rm dom}}\nolimits}
\def\hom{\mathop{\mathscr{H}}\nolimits}
\def\zero{\mathop{\mathscr{Z}}\nolimits}

\begin{document}
\vspace{5mm}

\date{\today}

\title[On derivations]{On derivations with respect to finite sets of smooth functions}

\author[R.\ Gr\"unwald]{Rich\'ard Gr\"unwald}
\author[Zs. P\'ales]{Zsolt P\'ales}
\address{Institute of Mathematics, University of Debrecen, 
H-4002 Debrecen, Pf. 400, Hungary}
\email{richard.grunwald96@gmail.com, pales@science.unideb.hu}

\thanks{The research of the first author was supported by the \'UNKP-17-1 New National Excellence Program of 
the Ministry of Human Capacities. The research of the second author was supported by the Hungarian Scientific Research Fund (OTKA) Grant
K-111651 and by the EFOP-3.6.1-16-2016-00022 project. This project is co-financed by the European Union and the European Social Fund.}

\begin{abstract}
The purpose of this paper is to show that functions that derivate the 
two-variable product function and one of the exponential, trigonometric or hyperbolic functions are also 
standard derivations. The more general problem considered is to describe finite sets of differentiable 
functions such that derivations with respect to this set are automatically standard derivations.
\end{abstract}

\maketitle

\section{Introduction}

Derivations are additive mappings of a ring into itself that possesses the so-called Leibniz Rule. More 
precisely, if $(R,+,\cdot)$ is a ring, then a function $d:R\to R$ is called a \emph{derivation} if, 
for all $x,y\in R$,
\begin{eqnarray}
  d(x+y)&=&d(x)+d(y), \label{ES} \\
   d(x\cdot y)&=&d(x)\cdot y+x\cdot d(y). \label{EP}
\end{eqnarray}
In other words, derivations behave similarly to the differentiation operator which is acting on 
differentiable real functions. A classical example of a derivation can be constructed on the ring $F[x]$ 
of polynomials over a field $F$ as follows: given $n\in\N$ and $a_0,\dots,a_n\in F$, define
\Eq{*}{
  d(a_nx^n+\dots+a_1x+a_0)=na_nx^{n-1}+\dots+a_1.
}
Then, an elementary calculation shows that $d:F[x]\to F[x]$ is a derivation.

Derivations are used in many branches of analysis and algebra. For instance, nonnegative information functions are constructed via real derivations (see  
Dar\'oczy--Maksa \cite{DarMak79}, Maksa \cite{Mak76}). Nonconstant functions that are convex with respect to families of power means are also obtained in terms of real derivations (see Maksa--P\'ales \cite{MakPal15}). Derivations are used to express the general solutions of certain functional equations (see Fechner--Gselmann \cite{FecGse12}, Gselmann \cite{Gse12}, \cite{Gse13a}, Halter-Koch \cite{Hal00}, \cite{Hal00b}, Jurkat \cite{Jur65}). Generalizations, such as higher-order derivations, bi-derivations and approximate or near-derivations were studied by Badora \cite{Bad06c}, Gselmann \cite{Gse14b}, Gselmann--P\'ales \cite{GsePal16}, and Maksa \cite{Mak81b}, \cite{Mak87c}.

In order to introduce the notion of a generalized derivation, first we define the classes of 
$n$-variable \emph{admissible functions} as follows:
\Eq{*}{
  \A_n:=\{f:\Omega\to\R \mid \Omega\subset\R^n \mbox{ is open, nonempty and }
          f\mbox{ is differentiable}\}
}
and
\Eq{*}{
  \A  := \bigcup_{n=1}^\infty \A_n.
}
The set $\Omega$ related to $f$ in the above definition will be called the domain of $f$ and denoted by 
$\dom_f$. It would make sense to introduce admissible functions as functions that map $\Omega$ into 
$\R^m$, but it is not necessary in this paper. We say that a function $d:\R\to\R$ is a \emph{derivation with 
respect to a subset $A\subseteq\A$} (shortly, \emph{d is a derivation with respect to $f_1$,\dots,$f_k$}, 
where $A=\{f_1,\dots,f_k\}$) if, for all $f\in A$,
\Eq{*}{
  d\big(f(x_1,\dots,x_n)\big)=\partial_1 f(x_1,\dots,x_n)d(x_1)+\dots+\partial_n f(x_1,\dots,x_n)d(x_n)
  \quad\,\big((x_1,\dots,x_n)\in\dom_f\big)
}
holds.
One can immediately see that a function $d:\R\to\R$ is a standard derivation if and 
only if it is a derivation with respect to $S_2$ and $P_2$, where 
\Eq{*}{
 S_2(x_1,x_2):=x_1+x_2 \qquad\mbox{and}\qquad P_2(x_1,x_2):=x_1x_2 \qquad((x_1,x_2)\in\R^2).
}

From now on we deal with functions that map $\R$ into $\R$. It is very simple to see some consequences of the Leibniz Rule. 
In the subsequent lemmas we describe the homogeneity properties of the solutions of the two functional 
equations \eq{S} and \eq{P}. We define the \emph{homogeneity set} and the \emph{set of zeros} of a 
function $d:\R\to\R$ by
\Eq{*}{
  \hom_d:=\{t\in\R\mid d(tx)=td(x) \mbox{ holds for all } x\in\R\} \qquad\mbox{and}\qquad
  \zero_d:=\{t\in\R\mid d(t)=0\},
}
respectively.

In the following three lemmas, we summarize the basic properties of various derivations.

\LEM{inv}{We have the following two assertions.
\begin{enumerate}[(i)]
 \item Let $\Omega_1$ and $\Omega_2$ be open subsets of $\R$, let $f:\Omega_1\to\Omega_2$ and $g:\Omega_2\to\R$ be differentiable functions. If $d:\R\to\R$ is a derivation with respect to $f$ and $g$, then $d$ is also a derivation with respect to the composition $g\circ f$.
 \item Let $\Omega_1$ and $\Omega_2$ be open subsets of $\R$, let $f:\Omega_1\to\Omega_2$ be a bijection such that $f'$ does not vanish on $\Omega_1$. If $d:\R\to\R$ is a derivation with respect to $f$, then $d$ is also a derivation with respect to the inverse function $f^{-1}:\Omega_2\to\Omega_1$.
\end{enumerate}}

\begin{proof}
By the assumptions of (i), for all $x\in\Omega_1$ and $y\in\Omega_2$, we have
\Eq{xy}{
	d(f(x))=f'(x)d(x)\qquad\mbox{and}\qquad d(g(y))=g'(y)d(y).
}
Therefore, with $y:=f(x)$, we get 
\Eq{*}{
  d((g\circ f)(x))=d(g(f(x)))=g'(f(x))d(f(x))=g'(f(x))f'(x)d(x)=(g\circ f)'(x)d(x),
}
which yields that $d$ is a derivation with respect to the function $g\circ f$.

By the assumption of (ii), for all $x\in\Omega_1$, we have the first equality in \eq{xy}.
Using the substitution $x=f^{-1}(y)$, this implies that
\Eq{*}{
	d(f^{-1}(y))=\frac{1}{f'(f^{-1}(y))}d(y)=\big(f^{-1}\big)'(y)d(y)\qquad(y\in \Omega_2).
} 
Thus, $d$ is a derivation with respect to the inverse function $f^{-1}$.
\end{proof}

\LEM{ADD}{Let $d:\R\to\R$ be a derivation with respect to $S_2$. Then $\hom_d$ is a 
subfield of $\R$, in particular $d$ is $\Q$-homogeneous, i.e., for all $r\in\Q$ and for all $x\in\R$,
\Eq{r}{
   d(rx)=rd(x).
}
Furthermore, $\zero_d$ is a vector space over $\hom_d$.
Additionally, if\, $\Omega\subseteq\R^n$ is a nonempty open set and $f,g:\Omega\to\R$ are differentiable 
functions such that $d$ is a derivation with respect to $f$ and $g$, then, for all $r,s\in\Q$, $d$ is a 
derivation with respect to $rf+sg$.}

\begin{proof}
One can easily see that $\hom_d$ is a subring of $\R$ with unit element $1$. Furthermore, for all 
$0\ne t\in\hom_d$ and $x\in\R$, $td\left(\frac{1}{t}x\right)=d\left(t\frac{1}{t}x\right)=d(x)$, which implies 
that $\frac{1}{t}\in\hom_d$. Therefore, $\hom_d$ is a subfield, indeed.

To justify that $\zero_d$ is a vector space over $\hom_d$, it suffices to observe that $\zero_d$ is a 
subgroup of $(\R,+)$, which is also closed by multiplications of elements of the field $\hom_d$.

Let $f,g:\Omega\to\R$ be differentiable functions such that $d$ is a derivation with respect to $f$ and $g$ 
and let $r,s\in\Q$. Then, using the additivity of $d$ and \eq{r} twice, we get
\Eq{*}{  
d((rf+sg)(x))
&=d(rf(x)+sg(x))=d(rf(x))+d(sg(x))\\
&=rd(f(x))+sd(g(x))=rf'(x)d(x)+sg'(x)d(x)=(rf+sg)'(x)d(x).
}
This shows that $d$ is a derivation with respect to $rf+sg$, indeed.
\end{proof}

\LEM{LR}{Let $d:\R\to\R$ be a derivation with respect to $P_2$. Then $\hom_d=\zero_d$ and 
$\hom_d\setminus\{0\}$ is a subgroup of $(\R,\cdot)$ containing the elements $-1,0$ and $1$. In particular, 
$d$ is odd, i.e., it is homogeneous with respect to $-1$. Furthermore, for all $r\in\Q$ and $x\in\R_+$,
\Eq{xr}{
  d(x^r)=rx^{r-1}d(x).
}
Additionally, if\, $\Omega\subseteq\R^n$ is a nonempty open set and $f,g:\Omega\to\R$ are differentiable 
functions such that $d$ is a derivation with respect to $f$ and $g$, then $d$ is a derivation with respect to 
$f\cdot g$. More generally, provided that $f$ and $g$ are positive on $\Omega$, for all $r,s\in\Q$, $d$ is 
also a derivation with respect to $f^r\cdot g^s$.}

\begin{proof}
Let $t\in \hom_d$. Then, for all $x\in\R$, by the Leibniz Rule, we get
\Eq{*}{
  td(x)=d(tx)=td(x)+xd(t).
}
Therefore, $d(t)=0$ must be valid, proving that $\hom_d\subseteq\zero_d$. The reversed inclusion follows 
similarly.

With the substitutions $x=y=0$, $x=y=1$ and $x=-y=1$ in equation \eq{P}, one can immediately see that 
$d(0)=d(1)=d(-1)=0$. This equation also shows that if $x,y\in \zero_d$, then $xy\in\zero_d$. Therefore, 
$\zero_d=\hom_d$ is a multiplicative subsemigroup of $\R$ containing $-1,0$ and $1$. For 
$0\neq x\in\zero_d$, we get that 
\Eq{*}{
  0=d(1)=d\Big(x\frac{1}{x}\Big)=xd\Big(\frac{1}{x}\Big)+\frac{1}{x}d(x),
}
hence $d\left(\frac{1}{x}\right)=0$. Therefore, $\hom_d\setminus\{0\}$ is a multiplicative subgroup of $\R$.

In this part, observe that with substitution $x:=e^u$, $y:=e^v$, the Leibniz Rule yields that 
the function $a:\R\to\R$ defined by
\Eq{*}{
  a(u):=e^{-u}d(e^u) \qquad(u\in\R)
}
is a derivation with respect to $S_2$. Therefore, by \lem{ADD}, $a$ is $\Q$-homogeneous. Hence
\Eq{*}{
e^{-ru}d(e^{ru})=re^{-u}d(e^u) \qquad(u\in\R,\,r\in\Q).
}
Substituting $u:=\ln x$ (where $x>0$) into the above identity, it follows that \eq{xr} holds.

Finally, we prove the last assertion of the lemma. Let $f,g:\Omega\to\R$ be positive differentiable functions 
such that $d$ is a derivation with respect to $f$ and $g$ and let $r,s\in\Q$. Then, using the Leibniz Rule
for $d$ and \eq{xr} twice, we get
\Eq{*}{  
d((f^r\cdot g^s)(x))
&=d(f^r(x)\cdot g^s(x))=g^s(x)d(f^r(x))+f^r(x)d(g^s(x))\\
&=rg^s(x)f^{r-1}(x)d(f(x))+sf^r(x)g^{s-1}(x)d(g(x)) \\
&=rg^s(x)f^{r-1}(x)f'(x)d(x)+sf^r(x)g^{s-1}(x)g'(x)d(x) \\
&=(rf^{r-1}f'\cdot g^s+f^r\cdot sg^{s-1}g')(x)d(x)=(f^rg^s)'(x)d(x).
}
This shows that $d$ is a derivation with respect to $f^r\cdot g^s$. If $r=s=1$, then the above 
argument can be applied not only for positive $f$ and $g$, hence we obtain that $d$ is a derivation with 
respect to the product function $f\cdot g$.
\end{proof}

The following result summarizes the most basic properties of real derivations. 

\THM{K}{(Kuczma \cite[pp.\ 346-352]{Kuc85}.) If $d:\R\to\R$ is a derivation, then $d(x)=0$ for every 
algebraic number $x\in\R$. On the other hand, for every non-algebraic number $t\in\R$, there exists a 
derivation $d:\R\to\R$ such that $d(t)\neq0$. If $d:\R\to\R$ is a derivation which is upper bounded on a set of 
positive measure, then $d$ is identically zero.}

A basic question that have been dealt with in several papers is to find conditions for additive functions 
which imply that this function is a standard derivation.

The following result is due to Nishiyama and Horinouchi \cite{NisHor68} (cf.\ Boros--Gselmann
\cite{BorGse10}).

\THM{NH}{Let $d:\R\to\R$ be a derivation with respect to $S_2$ and let $r\in\Q\setminus\{0,1\}$. Assume that 
$d$ satisfies the equation \eq{xr} for all $x\in\R_+$. Then $d$ is a standard derivation.}

The particular cases $r=-1$ and $r=2$ have been discovered by Kurepa \cite{Kur64}, \cite{Kur65}.
Results in the same spirit have also been established by Boros and Erdei \cite{BorErd05}.

\THM{BE}{Let $d:\R\to\R$ be a derivation with respect to $S_2$. Assume that $d$ satisfies the equation 
\Eq{*}{
  d(\sqrt{1-x^2})=-\frac{x}{\sqrt{1-x^2}}d(x)
}
for all $x\in]-1,1[$. Then $d$ is a standard derivation.}

The following theorem was proved by Gyula Maksa \cite{Mak13a}.

\THM{Mak}{Let $d:\R\to\R$ be a derivation with respect to $S_2$ and to one of the 
exponential, trigonometric or hyperbolic functions. Then $d$ is a standard derivation.}

Motivated by the above results, the purpose of this paper is to show that functions that derivate the 
two-variable product function and one of the exponential, trigonometric or hyperbolic functions are also 
standard derivations. The more general problem considered is to describe finite subsets of $\A$ 
such that derivations with respect to this set are automatically standard derivations.

\section{Main results}

The main assumption in this section is that $d:\R\to\R$ is a derivation with respect to $P_2$. Then, under 
various cicumstances, we prove that $d$ must be a standard derivation.

\subsection{The basic lemma}

The key result for our approach will be the Lemma formulated below. First, for $t\in\R\setminus\{0,-1\}$, 
define the set $H_t$ by
\Eq{*}{
  H_t:=\Big\{t,\frac1t,-1-\frac1t,-\frac{t}{1+t},-\frac{1}{1+t},-1-t\Big\}.
}
Observe that $H_t$ is a set which is invariant with respect to the mappings $\rho(s):=1/s$ and $\sigma(s):=-1-s$.
Furthermore, this set contains exactly six elements, unless $t\in\{1,-2,-\frac12\}$. In the latter case 
$H_t$ contains exactly three elements.

\Lem{BL}{Let $d:\R\to\R$ be a derivation with respect to $P_2$ and $U\subseteq\R$ be a set such that, for all 
$t\in\R\setminus\{0,-1\}$, the set $H_t$ intersects $U$. Assume that 
\Eq{d1}{
  d(u)=d(u+1) \qquad (u\in U).
}
Then $d$ is a standard derivation.}

\begin{proof}
We are going to show that the conditions of the lemma imply that if $x+y+z=0$, then $d(x)+d(y)+d(z)=0$ (where 
$x,y,z\in\R$). If $xyz=0$, then this statement is the consequence of the oddness of $d$. Therefore, 
we may assume that $xyz\neq0$. Using the condition $x+y+z=0$, observe that
\Eq{*}{
  \Big\{\frac{x}{y},\frac{y}{x},\frac{z}{x},\frac{x}{z},\frac{y}{z},\frac{z}{y}\Big\}=H_{\frac{x}{y}}.
}
Therefore, by our assumption, one of the above elements belongs to $U$. Due to the symmetry, we may assume 
that $u:=\frac{x}{y}\in U$. Then, by \eq{d1} and the Leibniz Rule, we obtain
\Eq{*}{
	0=d(u)-d(u+1)
	&=d\Big(\frac{x}{y}\Big)-d\Big(-\frac{z}y\Big) \\
	&=\frac{d(x)y-xd(y)}{y^2}+\frac{d(z)y-zd(y)}{y^2}
	=\frac{d(x)+d(y)+d(z)}{y}.
}
This implies that $d(x)+d(y)+d(z)=0$. In particular, if $x,y\in\R$ are arbitrary and $z=-x-y$, then we obtain 
$d(x)+d(y)+d(z)=0$, which implies that $d(x)+d(y)-d(x+y)=0$ and completes the proof of the additivity of 
$d$. 
Thus, $d$ is a standard derivation.
\end{proof}

\Cor{R}{Let $d:\R\to\R$ be a derivation with respect to $P_2$ and let $R\subseteq\,\, ]0,1[$ such that one of the following conditions holds:
\begin{enumerate}[(i)]
	\item $1\notin R+R$;
	\item $1\notin (\Z+R)(\Z+R)$;
	\item $-1\notin (\Z+R)(\Z+R)$. 
\end{enumerate}
Assume that
\Eq{*}{
	d(u)=d(u+1) \qquad (u\in\R\setminus(\Z+R)).
}
Then $d$ is a standard derivation.
} 

\begin{proof}In view of \lem{BL}, it suffices to show that, for all $t\in\R\setminus\{0,-1\}$, the set $H_t$ intersects $U:=\R\setminus(\Z+R)$. If condition (i) holds, then we show that $H_t\supset\{t,-1-t\}\cap U\neq\emptyset$. Indirectly suppose that $t\notin U$ and $-1-t\notin U$. Then $t\in\Z+R$ and $-1-t\in \Z+R$, hence $-1\in \Z+(R+R)$, which implies that $1\in R+R$ contradicting (i). Assume that condition (ii) holds. Then we prove that $H_t\supset\{t,\frac1t\}\cap U\neq\emptyset$. Similarly as before, if the previous intersection were empty, then $t\in\Z+R$ and $\frac1t\in\Z+R$, hence $1\in(\Z+R)(\Z+R)$ would be valid, which contradicts (ii). Finally, suppose that condition (iii) holds. In this case we show that $H_t\supset\{t,-1-\frac1t\}\cap U\neq\emptyset$. If this were not valid, then we would obtain that $t\in\Z+R$ and $-1-\frac1t\in\Z+R$. The last relation is equivalent to $-\frac1t\in\Z+R$. Thus, we would get $-1\in(\Z+R)(\Z+R)$, which is impossible by condition (iii). 
\end{proof}

\Cor{per}{
	Let $d:\R\to\R$ be a derivation with respect to $P_2$ and let further $0\leq 
	r<1$ be a constant. Assume that 
\Eq{*}{
	d(u)=d(u+1) \qquad (u\in\R\setminus(\Z+r)).
}
Then $d$ is a standard derivation.
}

\begin{proof}  If $0<r\neq\frac12$, then $1\notin \{r\}+\{r\}$, hence condition (i) of the previous corollary holds with $R:=\{r\}$, which yields the statement in this case. If $r=\frac12$, then $1\notin (\Z+\{r\})(\Z+\{r\})$ because otherwise, for some $n,k\in\Z$, we have that  $1=\big(n+\frac12\big)\big(k+\frac12\big)$, which is impossible. Thus, condition (ii) of the previous corollary holds with $R:=\{r\}$ again. Finally, suppose that $r=0$. If $t\in\R\setminus\Z$, then $t\in(\R\setminus\Z)\cap H_t$. If $t\in\Z$, then $\frac1t\in (\R\setminus\Z)\cap H_t$, unless $t=1$. In this remaining case $-\frac{1}{1+t}=-\frac12\in (\R\setminus\Z)\cap H_t$, which proves that $H_t$ intersects $\R\setminus\Z$ for all $t\neq0,-1$. Thus, the statement follows from  \lem{BL}.
\end{proof}

\Cor{int}{
	Let $d:\R\to\R$ be a derivation with respect to $P_2$. Denote by $I$ one of the following 6 intervals: 
	$I_1:=]-\infty,-2]$, $I_2:=[-2,-1[$, $I_3:=]-1-\frac12]$, $I_4:=[-\frac12,0[$, $I_5:=]0,1]$, and $I_6:=[1,\infty[$. Assume that 
\Eq{*}{
	d(u)=d(u+1) \qquad (u\in I).
}	
Then $d$ is a standard derivation.
}

\begin{proof} In view of \lem{BL}, it suffices to show that, for all $t\in\R\setminus\{0,-1\}$, the set $H_t$ 
intersects $U:=I$. Define the maps $\sigma:\R\to\R$ and $\rho:\R\setminus\{0\}\to\R\setminus\{0\}$ by 
$\sigma(t)=-1-t$ and $\rho(t)=1/t$, respectively. One can easily see that
\Eq{*}{
  &\sigma(I_1)=I_6,\quad && \sigma(I_2)=I_5,\quad && \sigma(I_3)=I_4,\quad && 
  \sigma(I_4)=I_3,\quad && \sigma(I_5)=I_2,\quad && \sigma(I_6)=I_1, \\
  &\rho(I_1)=I_4,\quad && \rho(I_2)=I_3,\quad && \rho(I_3)=I_2,\quad && 
  \rho(I_4)=I_1,\quad && \rho(I_5)=I_6,\quad && \rho(I_6)=I_5.
}
If $t\in\R\setminus\{-1,0\}$, then it belongs to one of the intervals $I_i$ ($i\in\{1,\dots,6\}$). Then one of the elements 
$t$, $\rho(t)=1/t$, $\sigma(\rho(t))=-1-1/t$, $\rho(\sigma(\rho(t)))=-t/(1+t)$, $\rho(\sigma(t))=-1/(1+t)$, and $\sigma(t)=-1-t$ will belong to $I$, proving that $I$ intersects $H_t$.
\end{proof}

\subsection{Periodic functions}

In the following result we deal with functions that are derivations with respect to $P_2$ and derivate a 
periodic or antiperiodic function. Given two constants $0\leq q<p$, we define the set $p\Z+q$ by
\Eq{*}{
  p\Z+q:=\{pk+q\mid k\in\Z\}.
}
We call a function $f:\R\setminus(p\Z+q)\to\R$ \emph{$p$-periodic} if, for all 
$x\in\R\setminus(p\Z+q)$,
\Eq{per}{
  f(x+p)=f(x).
}
We say that $f$ is \emph{$p$-antiperiodic} if, for all $x\in\R\setminus(p\Z+q)$,
\Eq{aper}{
  f(x+p)=-f(x).
}

\Thm{per}{
Let $d:\R\to\R$ be a derivation with respect to $P_2$, let $0\leq q<p$ be constants and let 
$f:\R\setminus(p\Z+q)\to\R$ be a differentiable $p$-periodic or 
$p$-antiperiodic function such that $f'(x)\neq0$ for all $x\in\R\setminus(p\Z+q)$. If  
\Eq{df}{
  d(f(x))=f'(x)d(x)
}
holds for all $x\in\R\setminus(p\Z+q)$, then $d$ is a standard derivation.
}

\begin{proof} In view of \lem{LR}, $d$ is odd and we have $d(0)=d(1)=d(-1)=0$.

Assume first that $f$ is $p$-periodic. Then also its derivative is $p$-periodic. Therefore, 
using \eq{df} twice, for all $x\in\R\setminus(p\Z+q)$, we get
\Eq{*}{
  f'(x)d(x)=d(f(x))=d(f(x+p))=f'(x+p)d(x+p)=f'(x)d(x+p).
}
Hence, for all $x\in\R\setminus(p\Z+q)$,
\Eq{dp}{
  d(x)=d(x+p),
}
i.e., $d$ is $p$-periodic.

If $f$ is $p$-antiperiodic, then $f'$ is also $p$-antiperiodic, therefore, using the oddness of $d$, we 
similarly get that
\Eq{*}{
  f'(x)d(x)=d(f(x))=d(-f(x+p))=-d(f(x+p))=-f'(x+p)d(x+p)=f'(x)d(x+p).
}
Thus, $d$ is $p$-periodic in this case, too.

If $0<q<p$, then $0\notin p\Z+q$, hence \eq{dp} with $x=0$ implies that $d(p)=0$. We show that $d(p)=0$ is 
also valid in the case $q=0$. In this case \eq{dp} holds for all $x\in\R\setminus p\Z$. In particular, it 
is valid for $x=-\frac{p}{2}$, hence
\Eq{*}{
   d\Big(-\frac{p}{2}\Big)=d\Big(\frac{p}{2}\Big),
}
which implies that $d\big(\frac{p}{2}\big)=0$, whereby, using \eq{dp}, 
$d\big((2k+1)\frac{p}{2}\big)=0$ holds for all $k\in\Z$. By \lem{LR}, the zeros of $d$ form a multiplicative 
subgroup in $\R$, hence $d\big(\frac{2k+1}{2\ell+1}\big)=0$ is valid for 
all $k,\ell\in\Z$. Now using the previous equation and substituting $x=-\frac{p}{4}$ into \eq{dp}, we obtain
\Eq{*}{
  -d\Big(\frac{p}{4}\Big)=d\Big(\frac{-p}{4}\Big)=d\Big(\frac{3p}{4}\Big)=3d\Big(\frac{p}{4}\Big),
}
hence $d\left(\frac{p}{4}\right)=0$. Thus, 
\Eq{*}{
	d(2)=d\left(\frac{\frac{p}{2}}{\frac{p}{4}}\right)=0,
} 
whence 
\Eq{*}{
  d(p)=d\Big(2\frac{p}{2}\Big)=2d\Big(\frac{p}{2}\Big)+\frac{p}{2}d(2)=0.
}
Then, in each of the two cases, we have that $d$ is $p$-homogeneous. 

Let $y\in\R\setminus\big(\Z+\frac{q}{p}\big)$ be arbitrary. Then $py\in \R\setminus(p\Z+q)$, therefore, by 
\eq{dp}, we obtain that
\Eq{*}{
  d(y)=\frac1pd(py)=\frac1pd(py+p)=d(y+1)
}
for all $y\in\R\setminus\big(\Z+\frac{q}{p}\big)$. Using \cor{per}, it follows that $d$ is a standard 
derivation.
\end{proof}

In the following consequence of \thm{per} we obtain several equivalent conditions in terms of the 
derivation of trigonometric functions.

\Cor{tri}{Let $d:\R\to\R$ be a derivation with respect to $P_2$. Then the following statements are equivalent:
\begin{enumerate}[(i)]
	\item $d$ is a derivation with respect to the sine function;
	\item $d$ is a derivation with respect to the cosine function;
	\item $d$ is a derivation with respect to the tangent function;
	\item $d$ is a derivation with respect to the cotangent function.
\end{enumerate}
Furthermore, in each of the above cases, $d$ is a standard derivation.}

\begin{proof} 
Let $d:\R\to\R$ be a derivation with respect to $P_2$ and assume that condition (i) holds, that 
is,
\Eq{sin}{
   d(\sin(x))=\cos(x)d(x) \qquad(x\in\R).
}
In order to apply \thm{per}, define $f$ to be the restriction of sine to the set 
$\R\setminus(\pi\Z+\frac{\pi}{2})$. Then $f$ is $\pi$-antiperiodic and $f'$ does not vanish on 
$\R\setminus(\pi\Z+\frac{\pi}{2})$. The equation \eq{sin} yields that \eq{df} also holds for 
$x\in\R\setminus(\pi\Z+\frac{\pi}{2})$. Therefore, by \thm{per}, it follows that $d$ is a standard derivation.
Substituting $x=\pi$ into \eq{sin}, we get that $d(\pi)=0$. By the additivity of $d$, \lem{ADD} implies that 
$d\big(\frac\pi2\big)=0$. Therefore, using \eq{sin}, for all $x\in\R$, we obtain
\Eq{*}{
  d(\cos(x))
  =d\Big(\sin\Big(x+\frac\pi2\Big)\Big)
  &=\cos\Big(x+\frac\pi2\Big)d\Big(x+\frac\pi2\Big)\\
  &=-\sin(x)\Big(d(x)+d\Big(\frac\pi2\Big)\Big)=-\sin(x)d(x),
}
which proves that condition (ii) also holds. Based on \lem{LR}, then $d$ is 
a derivation with respect to the tangent and cotangent functions, i.e., conditions (iii) and (iv) also follow from (i).

In the second part of the proof, suppose that $d:\R\to\R$ is a derivation with respect to $P_2$ and condition (ii) holds, that is,
\Eq{cos}{
   d(\cos(x))=-\sin(x)d(x) \qquad(x\in\R).
}
Defining the function $f$ as the restriction of cosine to the set $\R\setminus(\pi\Z)$, \thm{per} implies 
that $d$ is a standard derivation. Then, following a similar train of thought as above, we can 
get that conditions (i), (iii), and (iv) are also valid. 

In the third part, assume that $d:\R\to\R$ is a derivation with respect to $P_2$ and condition (iii) holds, that is,
\Eq{tan}{
   d(\tan(x))=\frac{d(x)}{\cos^2(x)} \qquad(x\in\R\setminus(\pi\Z+\tfrac{\pi}{2})).
}
Applying \thm{per} for $f:=\tan$, it immediately follows that $d$ is a standard derivation.
Substituting $x=\pi$ into \eq{tan}, we get that $d(\pi)=0$. By the additivity of $d$, \lem{ADD} implies that 
$d\big(\frac\pi2\big)=0$. Therefore, using that $d$ is a standard derivation and that \eq{tan} 
holds, for all $x\in\R\setminus(2\pi\Z+\pi)$, we obtain
\Eq{*}{
  d(\sin(x))
  =d\Big(\frac{2\tan(x/2)}{1+\tan^2(x/2)}\Big)
  &=\frac{2(1-\tan^2(x/2))}{(1+\tan^2(x/2))^2}d(\tan(x/2))\\
  &=\frac{2(1-\tan^2(x/2))}{(1+\tan^2(x/2))^2}\cdot\frac{d(x/2)}{\cos^2(x/2)}
  =\cos(x)d(x),
}
which proves that condition (i) is satisfied.

In the fourth part, suppose that $d:\R\to\R$ is a derivation with respect to $P_2$ and condition (iv) holds, that is,
\Eq{cot}{
	d(\cot(x))=-\frac{d(x)}{\sin^2(x)} \qquad(x\in\R\setminus(\pi\Z)).
}
Using \thm{per} for $f:=\cot$, it immediately follows that $d$ is a standard derivation. Then, following a 
similar train of thought as in the third part, we can get that conditions (i) is also valid. 
\end{proof}

\subsection{Exponential and hyperbolic functions}

\Thm{hyp}{Let $d:\R\to\R$ be a derivation with respect to $P_2$. Assume that 
$d(2)=0$. Then the following statements are equivalent:
\begin{enumerate}[(i)]
	\item $d$ is a derivation with respect to the sine hyperbolic function;
	\item $d$ is a derivation with respect to the cosine hyperbolic function;
	\item $d$ is a derivation with respect to the tangent hyperbolic function;
	\item $d$ is a derivation with respect to the cotangent hyperbolic function;
	\item $d$ is a derivation with respect to the exponential function.
\end{enumerate}
Furthermore, in each of the above cases, $d$ is a standard derivation.}

\begin{proof} 
Let $d:\R\to\R$ be a derivation with respect to $P_2$ and assume that condition (i) holds, that is,
\Eq{*}{
 d(\sinh(x))=\cosh(x)d(x) \qquad(x\in\R).
}
The property $d(2)=0$ implies that $d$ is $2$-homogeneous. 
Using the Leibniz Rule and that $\sinh(2x)=2\sinh(x)\cosh(x)$ and $\cosh(2x)=\cosh^2(x)+\sinh^2(x)$ hold for 
all $x\in\R$, we obtain
\Eq{*}{
  2(\sinh(x)d(\cosh(x))+\cosh^2(x)d(x))
  &=2(\sinh(x)d(\cosh(x))+\cosh(x)d(\sinh(x)))\\ 
  &=d(2\sinh(x)\cosh(x))=d(\sinh(2x))\\&=\cosh(2x)d(2x)=2(\cosh^2(x)+\sinh^2(x))d(x).
} 
Hence
\Eq{*}{
 \sinh(x)d(\cosh(x))=\sinh^2(x)d(x)
}
holds. This equality, for $x\neq0$, simplifies to
\Eq{*}{
 d(\cosh(x))=\sinh(x)d(x),
}
which is also valid for $x=0$. Therefore, condition (i) implies condition (ii). Based on \lem{LR}, then $d$ is 
a derivation with respect to the tangent and cotangent hyperbolic functions, i.e., conditions (iii) and (iv) 
also follow from (i). Since $\exp=\sinh+\cosh$, thus $d$ is also a derivation with respect to the exponential 
function, i.e., (v) holds.

In the second part, we first show that condition (ii) implies that $d$ is a standard derivation and 
then that conditions (i), (iii) and (iv) also consequences of (ii). Using the assumptions and some well-known 
identities, we get
\Eq{*}{
  d(2\cosh^2(x)-1)&=d(\cosh(2x))=\sinh(2x)d(2x) \\& =4\sinh(x)\cosh(x)d(x)=4\cosh(x)d(\cosh(x))=d(2\cosh^2(x))
}
holds for all $x\in\R$. Substituting $u:=2\cosh^2(x)-1$, we obtain
\Eq{*}{
 d(u)=d(u+1) \qquad (u\geq1),
}
hence, in view of \cor{int}, $d$ is a standard derivation. To justify that condition (ii) implies (i), we use 
\lem{LR} and that $d$ is a derivation. We have that
\Eq{*}{
 d(\sinh(x))=d\Big(\sqrt{\cosh^2(x)-1}\Big)=\frac{1}{2}\frac{d(\cosh^2(x)-1)}{\sqrt{\cosh^2(x)-1}}
 =\frac{1}{2}\frac{2\cosh(x)\sinh(x)d(x)}{\sinh(x)}=\cosh(x)d(x)
}
is valid for all $x\geq0$. By the oddness of $d$, the above equality also holds for $x<0$ and therefore $d$ 
is a derivation with respect to the tangent and cotangent hyperbolic functions, i.e., conditions (iii) 
and (iv) also follow from (ii). Since $\exp=\sinh+\cosh$, thus $d$ is also a derivation with respect to the 
exponential function.

In the third part, observe that, in view of \lem{LR}, (iii) and (iv) are equivalent to each other. Thus, 
it suffices to show that condition (iv) implies that $d$ is a standard derivation and then that conditions 
(i) and (ii) are also consequences of (iv). By the assumption and some well-known identities for hyperbolic 
functions, we obtain
\Eq{*}{
  d\Big(\frac{\coth^2(x)+1}{2\coth(x)}\Big) 
  =d(\coth(2x))
  &=-\frac{1}{\sinh^2(2x)}d(2x)
  =-\frac{1}{4\sinh^2(x)\cosh^2(x)}2d(x) \\& 
  =\frac{1}{2\cosh^2(x)}d(\coth(x))
  =\frac{\coth^2(x)-1}{2\coth^2(x)}d(\coth(x))
}
holds for all $x\in\R$. Substituting $t:=\coth(x)$, we get
\Eq{*}{
  d\Big(\frac{t^2+1}{2t}\Big)=\frac{t^2-1}{2t^2}d(t) \qquad (|t|>1).
}
Using the Leibniz Rule again, we obtain
\Eq{*}{
  \frac{(t^2+1)d(t)-td(t^2+1)}{t^2}=\frac{t^2-1}{t^2}d(t) \qquad (|t|>1).
}
Therefore, after reduction, we have that
\Eq{*}{
	d(t^2)=d(t^2+1)	 \qquad (|t|>1),
}
hence with $u:=t^2$,
\Eq{*}{
	d(u)=d(u+1)  \qquad (u>1).
}
The above equality also holds for $u=1$, thus $d$ is a standard derivaton by \cor{int}. Using this, we 
show that condition (iv) implies (i). Indeed, 
\Eq{*}{
	d(\sinh(x))=d\left(\frac{1}{\sqrt{\coth^2(x)-1}}\right)
	&=-\frac{1}{2}\frac{2\coth(x)}{\sqrt{(\coth^2(x)-1)^3}}\frac{-1}{\sinh^2(x)}d(x)=\cosh(x)d(x)
}
is valid for all $x>0$. By the oddness of $d$, the above equality also holds for $x\leq0$. A similar 
computation yields that (ii) is also a consequence of property (iv). Since $\exp=\sinh+\cosh$, thus it 
follows that $d$ is also a derivation with respect to the exponential function.

To complete the proof, assume that condition (v) holds, i.e.,
\Eq{E}{
	d(e^x)=e^xd(x)\qquad(x\in\R).
}
Then, using the Leibniz Rule and \eq{E} three times, for $x,y\in\R$, we obtain
\Eq{*}{
	d(x+y)&=e^{-x-y}d(e^{x+y})
	=e^{-x-y}d(e^{x}e^{y}) 
	=e^{-x-y}\big(e^yd(e^x)+e^xd(e^y)\big)\\
	&=e^{-x-y}\big(e^ye^xd(x)+e^xe^yd(y)\big)=d(x)+d(y).
}
Therefore, $d$ is a derivation with respect to $S_2$, thus it is a 
standard derivation, indeed. Since $\sinh(x)=\frac12(e^x-e^{-x})$, therefore $d$ is a derivation with respect 
to the sine hyperbolic function, whence (i)-(iv) follow immeadiately.
\end{proof}

\Rem{exp}{In the above theorem, for those implications, where condition (v) is supposed, it is not necessary to assume that $d(2)=0$.
The question whether the above theorem can be proved without the condition $d(2)=0$ remains open.}


\begin{thebibliography}{10}

\bibitem{Bad06c}
R.~Badora, On approximate derivations, \emph{Math. Inequal. Appl.}, \textbf{9} (2006), 167--173. 

\bibitem{BorErd05}
Z.~Boros and P.~Erdei, A conditional equation for additive functions, \emph{Aequationes Math.}, \textbf{70} (2005), 309--313. 

\bibitem{BorGse10}
Z.~Boros and E.~Gselmann, Hyers--{U}lam stability of derivations and linear functions, \emph{Aequationes Math.}, \textbf{80} (2010), 13--25.

\bibitem{DarMak79}
Z.~Dar\'oczy and Gy. Maksa, Nonnegative information functions, \emph{Analytic function methods in probability theory (Proc. Colloq. Methods of
  Complex Anal. in the Theory of Probab. and Statist., Lajos Kossuth Univ. Debrecen, Debrecen, 1977)}, (North-Holland, Amsterdam, 1979), p.~67--78.

\bibitem{FecGse12}
W.~Fechner and E.~Gselmann, General and alien solutions of a functional equation and of a functional inequality, \emph{Publ. Math. Debrecen}, \textbf{80}
  (2012), 143--154.

\bibitem{Gse12}
E.~Gselmann, Notes on the characterization of derivations, \emph{Acta Sci. Math. (Szeged)}, \textbf{78} (2012), 137--145.

\bibitem{Gse13a}
E.~Gselmann, Derivations and linear functions along rational functions, \emph{Monatsh. Math.}, \textbf{169} (2013), 355--370.

\bibitem{Gse14b}
E.~Gselmann, Approximate derivations of order {$n$}, \emph{Acta Math. Hungar.}, \textbf{144} (2014), 217--226. 

\bibitem{GsePal16}
E.~Gselmann and Zs. P\'ales, Additive solvability and linear independence of the solutions of a system of functional equations \emph{Acta Sci. Math.
  (Szeged)}, \textbf{82} (2016), 101--110. 

\bibitem{Hal00b}
F.~Halter-Koch, A characterization of derivations by functional equations, \emph{Math. Pannon.}, \textbf{11} (2000), 187--190. 

\bibitem{Hal00}
F.~Halter-Koch, Characterization of field homomorphisms and derivations by functional equations, \emph{Aequationes Math.}, \textbf{59} (2000), 298--305. 

\bibitem{Jur65}
W.~B. Jurkat, On {C}auchy's functional equation, \emph{Proc. Amer. Math. Soc.}, \textbf{16} (1965), 683--686. 

\bibitem{Kuc85}
M.~Kuczma, \emph{{An {I}ntroduction to the {T}heory of {F}unctional {E}quations and {I}nequalities}}, {Prace Naukowe Uniwersytetu \'Sl\c{a}skiego w Katowicach}, vol. 489, Pa\'nstwowe Wydawnictwo Naukowe -- Uniwersytet \'Sl\c{a}ski, (Warszawa--Krak\'ow--Katowice, 1985), 2nd edn. (ed. by A. Gil\'anyi),
  (Birkh\"auser, Basel, 2009). 

\bibitem{Kur64}
S.~Kurepa, The {C}auchy functional equation and scalar product in vector spaces, \emph{Glasnik Mat.-Fiz. Astronom. Dru\v{s}tvo Mat. Fiz. Hrvatske Ser. II}, 
  \textbf{19} (1964), 23--36. 

\bibitem{Kur65}
S.~Kurepa, Remarks on the {C}auchy functional equation \emph{Publ. Inst. Math. (Beograd) (N.S.)},  \textbf{5 (19)} (1965), 85--88. 

\bibitem{Mak76}
Gy. Maksa, \emph{Devi\'aci\'ok \'es differenci\'ak ({D}eviations and
  {D}ifferences)}, (in {H}ungarian), Phd thesis, (Lajos Kossuth University, Debrecen, Hungary, 1976).

\bibitem{Mak81b}
Gy. Maksa, On near derivations, \emph{Proc. Amer. Math. Soc.}, \textbf{81} (1981), 406--408. 

\bibitem{Mak87c}
Gy. Maksa, On the trace of symmetric bi-derivations, \emph{ C. R. Math. Rep. Acad. Sci. Canada}, \textbf{9} (1987), 303--307.

\bibitem{Mak13a}
Gy. Maksa, On additive functions which differentiate elementary functions in some sense, \emph{Ann. Univ. Sci. Budapest. Sect. Comput.}, \textbf{41} (2013),
  125--136. 

\bibitem{MakPal15}
Gy. Maksa and Zs. P\'ales, Convexity with respect to families of means, \emph{Aequationes Math.}, \textbf{89} (2015), 161--167. 

\bibitem{NisHor68}
A.~Nishiyama and S.~Horinouchi, On a system of functional equations, \emph{Aequationes Math.}, \textbf{1} (1968), 1--5. 

\end{thebibliography}

\end{document}